\documentclass[runningheads]{llncs}
\usepackage{color}
\usepackage{amsfonts}
\usepackage{amsmath}
\usepackage{amssymb}
\usepackage{times}
\usepackage{cite}
\usepackage{hhline}
\usepackage{compactbib}
\usepackage{float}
\usepackage{calc}
\usepackage{enumerate}
\usepackage{subfig}
\usepackage{graphicx}

\renewcommand{\thefigure}{\arabic{figure}}
\renewcommand{\thetable}{\Roman{table}}

\usepackage[english]{babel}

\newtheorem{thm}{Theorem}[section]
\newtheorem{defn}[thm]{Definition}

\usepackage[numbers]{natbib}
\usepackage{tikz}
\usetikzlibrary{calc,patterns,decorations.pathmorphing,decorations.markings}
\usetikzlibrary{arrows,calc,shapes,decorations.pathreplacing}
\usetikzlibrary{decorations.markings}
\usetikzlibrary{trees,positioning,fit,shapes}
\usetikzlibrary{backgrounds}
\usetikzlibrary{patterns, positioning, arrows}

\makeatletter
\newcommand\setveclength[3]{
  \pgfpointdiff{\pgfpointanchor{#2}{center}}{\pgfpointanchor{#3}{center}}
  \pgfmathveclen{\pgf@x}{\pgf@y}
  \edef#1{\pgfmathresult}
}
\makeatother

\usepackage{amsmath}
\usepackage{bm}
\usepackage{amsfonts}
\usepackage{xfrac}
\usepackage{graphicx}

\begin{document}

\title{A topological analysis of cointegrated data: a Z24 Bridge case study}


\author{T. Gowdridge \and
E.J. Cross \and
N. Dervilis \and
K. Worden \\
\email{tgowdridge1@sheffield.ac.uk}}
\authorrunning{T. Gowdridge et al.}
\institute{Dynamics Research Group, Department of Mechanical Engineering, \\ University of Sheffield \\
Mappin Street, Sheffield S1 3JD, UK \\
}

\maketitle

\begin{abstract}
The paper studies the topological changes from before and after cointegration, for the natural frequencies of the Z24 Bridge. The second natural frequency is known to be nonlinear in temperature, and this will serve as the main focal point of this work. Cointegration is a method of normalising time series data with respect to one another - often strongly-correlated time series. Cointegration is used in this paper to remove effects from Environmental and Operational Variations, by cointegrating the first four natural frequencies for the Z24 Bridge data. The temperature effects on the natural frequency data are clearly visible within the data, and it is desirable, for the purposes of structural health monitoring, that these effects are removed. The univariate time series are embedded in higher-dimensional space, such that interesting topologies are formed. Topological data analysis is used to analyse the raw time series, and the cointegrated equivalents. A standard topological data analysis pipeline is enacted, where simplicial complexes are constructed from the embedded point clouds. Topological properties are then calculated from the simplicial complexes; such as the persistent homology. The persistent homology is then analysed, to determine the topological structure of all the time series.

\end{abstract}

\section {Introduction}

Cointegration has previously been used in Structural Health Monitoring (SHM) to remove the effects of Environmental and Operational Variations (EOVs) in collected data \cite{cross2011cointegration}. One application of cointegration was in removing such EOVs from the Z24 Bridge data set \cite{shi2016nonlinear, shi2018nonlinear, shi2018regime, shi2019cointegration}. When cointegrating data, the aim is to remove the long-term correlation between one or more time series. The effect on the topology of the time series has never been discussed before. This paper will explore the effects of cointegration and how the topology is affected, before and after.

The second natural frequency from the Z24 data set is the main focus for this paper, as this  presents strong nonlinear characteristics as a result of freezing temperatures. An assessment of the shortcomings of applying a linear cointegration scheme to a nonlinear problem will be explored. The linear natural frequencies are also briefly discussed.

Two cointegration schemes will be used in this paper: linear and nonlinear. The linear scheme assumes that $\omega_2$ is linearly correlated to the other natural frequencies; this is not actually true here, as there are nonlinear complications arising from freezing temperatures over certain measurement readings. To account for nonlinear effects, in the case of the nonlinear cointegration scheme, Gaussian Processes (GPs) are fitted to the natural frequency, these are then used in the calculations.

To form a topology from the univariate time series, a time-delay embedding is used \cite{takens1981detecting, packard1980geometry}. This paper extends on ideas of applying Topological Data Analysis (TDA) to the Z24 Bridge data, previously considered in  \cite{gowdridge2021z24}.

A standard TDA procedure is followed: simplicial complexes are constructed from some point cloud data, then the persistent homology is calculated. The persistent homology is used to compare different point clouds. 

The outline of this paper is as follows: Section 2, will give a brief walk-through of the background theory, pertinent to cointegration, TDA, and time-delay embeddings. Section 3 will then focus on analysis of the Z24 cointegrated data. Section 4 will provide a discussion of the results, and then the paper will conclude.

\section{Background Theory}

\subsection{Linear Cointegration}
Only a brief and qualitative introduction to cointegration will be given here, as cointegration is a previously well-covered topic in SHM. For greater detail on cointegration theory, the authors refer to other introductions \cite{cross2011cointegration, shi2019cointegration, dickey1979distribution}. Linear cointegration has previously been used to remove EOVs from the DAMASCOS data-set \cite{cross2011cointegration},  proving its use in novelty detection, among many other applications. 

Cointegration is a property of multiple nonstationary time series. If a selection of time series have this property, at least one linear combination of them exists that results in a stationary time series, referred to as the \textit{residual series}. This combination of time series that results in a stationary time series are encoded in the \textit{cointegrating vectors}. For a set of nonstationary time series $\{y_i(t)\}$, these are cointegrated iff,
\begin{equation}
    z_i(t) = \{\beta\}^{T}\{y_i(t)\},
\end{equation}
where $\{\beta\}$ is the cointegrating vector and $z_i(t)$ is a stationary time series.

Firstly, nonstationary time series are modelled with error-correction models, using a least-squares approach, shown below \cite{cross2011cointegration},

\begin{equation}
\label{eq:ecm}
    \nabla y_i = \rho y_{i-1} + \sum b_j \nabla y_{i-j} + \epsilon_i
\end{equation}
where $\nabla$ is the difference operator, and $\epsilon$ is the time-series residual, and should be a white-noise process.

For time series to be cointegrated, they must have the same degree of nonstationarity, this means that they must be differenced of the same order, to result in a stationary series. Once a nonstationary time series has been differenced, the Augmented Dickey-Fuller (ADF) test \cite{dickey1979distribution, dickey1981likelihood} must be employed to prove that the now-differenced time series can be considered stationary. The ADF test is a hypothesis test based on values of $\rho$ in the error-correction model shown in equation (\ref{eq:ecm}),
\begin{equation}
    t_p = \frac{\hat{\rho}}{\sigma_\rho}
\end{equation}
where $\hat{\rho}$ is the least squares estimate and $\sigma_\rho$ is the variance of the parameter. $t_p$ is compared against previously-calculated reference values \cite{dickey1979distribution, dickey1981likelihood}. The ADF test is essentially seeing how close $\rho$ is to 0. If $\rho=0$, the time series has a unit root, and is inherently nonstationary \cite{cross2011cointegration}. If $\rho$ is statistically close enough to 0, the test statistic is rejected, and the time series is considered nonstationary. The time series is then required to be differenced again. This process is repeated until the null hypothesis is accepted. The number of differences required for the time series to become stationary, is its order of integration.

To solve for the cointegrating vectors, the Johansen procedure \cite{cross2011cointegration, johansen1995likelihood} is followed here, this is essentially an eigenvalue problem; the cointegrating vectors are the eigenvectors. The cointegrating vector that results in the most stationary combination is the vector associated to the largest eigenvalue.

\subsection{Nonlinear Cointegration}
Linear cointegration has been developed for use on linear time series. When nonlinearities are present in a time series, no linear combination of the time series will result in a stationary residual \cite{cross2011approaches}; Thus, invalidating the assumption that the residual series is sampled from a white-noise process.

The general procedure to perform nonlinear cointegration begins similarly to linear cointegration. The ADF test should be implemented to ensure that the variables are integrated of the same order. The data sets should then be split into an training (containing no damage data) and test data sets used to train a GP regression model. The ADF test should then be applied to the trained GP to ensure that the residual series is of a lower order of integration, than when compared to linear cointegration \cite{shi2019cointegration, shi2016nonlinear, shi2018nonlinear}.

\subsection{Topological Data Analysis}

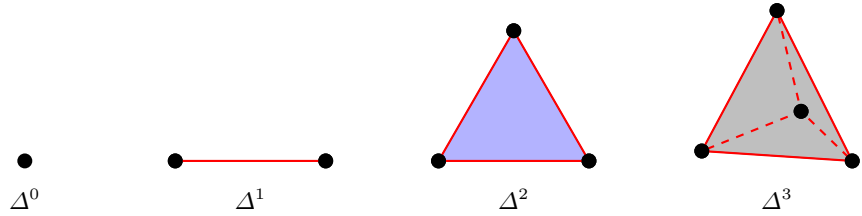
\begin{figure}
    \centering
\begin{tikzpicture}[ele/.style={fill=black,circle,minimum width=3pt,inner sep=2pt}]
\node[ele] (a1) at (0,0) {};
\node at (0,-0.5) {$\Delta^0$};

\begin{scope}[shift={(2,0)}]
    \node[ele] (a1) at (0,0) {};    
    \node[ele] (a2) at (2,0) {};
    \draw[-,thick,shorten <=2pt,shorten >=2pt, red] (a1.center) -- (a2.center);
    \node at (1,-0.5) {$\Delta^1$};
    
    \node[ele] (a1) at (0,0) {};    
    \node[ele] (a2) at (2,0) {};
\end{scope}

\begin{scope}[shift={(5.5,0)}]
    \node[ele] (a1) at (2,0) {};    
    \node[ele] (a2) at (0,0) {};
    \node[ele] (a3) at (1,1.732) {};

    \draw[-,thick,shorten <=2pt,shorten >=2pt, red] (a1.center) -- (a2.center);
    \draw[-,thick,shorten <=2pt,shorten >=2pt, red] (a1.center) -- (a3.center);
    \draw[-,thick,shorten <=2pt,shorten >=2pt, red] (a2.center) -- (a3.center);

    \begin{scope}[on background layer]
        \path [fill=blue!30, draw] (a1.center) to (a2.center) to (a3.center) to (a1.center);
    \end{scope} 

    \node at (1,-0.5) {$\Delta^2$};
    \node[ele] (a1) at (2,0) {};    
    \node[ele] (a2) at (0,0) {};
    \node[ele] (a3) at (1,1.732) {};
    
\end{scope}

\begin{scope}[shift={(9,0)}]
    \node[ele] (a1) at (2,0) {};    
    \node[ele] (a2) at (0,0.132) {};
    \node[ele] (a3) at (1.32,0.66) {};
    \node[ele] (a4) at (1,2) {};
    
    \draw[-,thick,shorten <=2pt,shorten >=2pt, red] (a1.center) -- (a2.center);
    \draw[-,thick,shorten <=2pt,shorten >=2pt, red] (a1.center) -- (a4.center);
    \draw[-,thick,shorten <=2pt,shorten >=2pt, red] (a2.center) -- (a4.center);

    \draw[-,thick,shorten <=2pt,shorten >=2pt, dashed, red] (a1.center) -- (a3.center);
    \draw[-,thick,shorten <=2pt,shorten >=2pt, dashed, red] (a2.center) -- (a3.center);
    \draw[-,thick,shorten <=2pt,shorten >=2pt, dashed, red] (a4.center) -- (a3.center);
    
    \begin{scope}[on background layer]
        \path [fill=lightgray, draw] (a1.center) to (a2.center) to (a4.center) to (a1.center);
    \end{scope} 
\node at (1,-0.5) {$\Delta^3$};
    \node[ele] (a1) at (2,0) {};    
    \node[ele] (a2) at (0,0.132) {};
    \node[ele] (a3) at (1.32,0.66) {};
    \node[ele] (a4) at (1,2) {};

\end{scope}
\end{tikzpicture}
    \caption{The first four simplices.}
    \label{fig:simplices}
\end{figure}

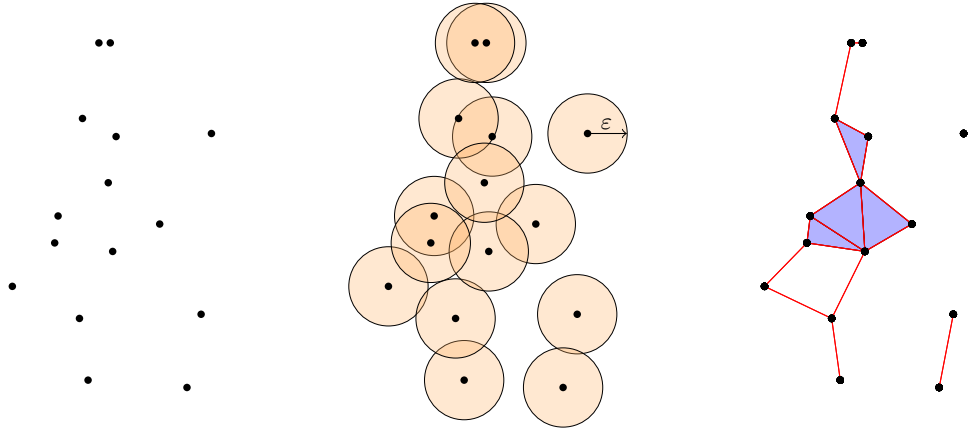
\begin{figure}
    \centering
    \begin{tikzpicture}[ele/.style={fill=black,circle,minimum width=2pt,inner sep=1pt}, scale = 1]

    \node[ele] (a) at (2.8209087137275, 1.18093727264829) {};
    \node[ele] (b) at (2.95833231728538,	3.58594837328201) {};
    \node[ele] (c) at (0.311441915093690,	1.55228774215654) {};
    \node[ele] (d) at (1.69060866126193,	3.54516924716448) {};
    \node[ele] (e) at (1.31799287927655,	0.304557156170160) {};
    \node[ele] (f) at (1.61347687509477,	4.79163906483583) {};
    \node[ele] (g) at (2.27043233468693,	2.38135190156370) {};
    \node[ele] (h) at (1.20381702801476,	1.12542610389225) {};
    \node[ele] (i) at (0.919304590945695,	2.48862221210784) {};
    \node[ele] (j) at (1.64512903092519,	2.01677357394698) {};
    \node[ele] (k) at (1.24383442716840,	3.78580579686650) {};
    \node[ele] (l) at (1.58646607060451,	2.93089871716034) {};
    \node[ele] (m) at (0.874518089942222,	2.13017381291781) {};
    \node[ele] (n) at (1.46214937616206,	4.79124963787566) {};
    \node[ele] (o) at (2.63312304723148,	0.206799557528458) {};

\begin{scope}[elem/.style={circle, draw = black, line width = 0.2pt, fill = orange!60, fill opacity=0.3, minimum width = 30pt}, shift={(5,0)}]

\node[elem] (a) at (2.8209087137275, 1.18093727264829) {};
\node[elem] (b) at (2.95833231728538,   3.58594837328201) {};
\node[elem] (c) at (0.311441915093690,  1.55228774215654) {};
\node[elem] (d) at (1.69060866126193,   3.54516924716448) {};
\node[elem] (e) at (1.31799287927655,   0.304557156170160) {};
\node[elem] (f) at (1.61347687509477,   4.79163906483583) {};
\node[elem] (g) at (2.27043233468693,   2.38135190156370) {};
\node[elem] (h) at (1.20381702801476,   1.12542610389225) {};
\node[elem] (i) at (0.919304590945695,  2.48862221210784) {};
\node[elem] (j) at (1.64512903092519,   2.01677357394698) {};
\node[elem] (k) at (1.24383442716840,   3.78580579686650) {};
\node[elem] (l) at (1.58646607060451,   2.93089871716034) {};
\node[elem] (m) at (0.874518089942222,  2.13017381291781) {};
\node[elem] (n) at (1.46214937616206,   4.79124963787566) {};
\node[elem] (o) at (2.63312304723148,   0.206799557528458) {};

\node[ele] (a1) at (2.8209087137275, 1.18093727264829) {};
\node[ele] (b1) at (2.95833231728538,    3.58594837328201) {};
\node[ele] (c1) at (0.311441915093690,   1.55228774215654) {};
\node[ele] (d1) at (1.69060866126193,    3.54516924716448) {};
\node[ele] (e1) at (1.31799287927655,    0.304557156170160) {};
\node[ele] (f1) at (1.61347687509477,    4.79163906483583) {};
\node[ele] (g1) at (2.27043233468693,    2.38135190156370) {};
\node[ele] (h1) at (1.20381702801476,    1.12542610389225) {};
\node[ele] (i1) at (0.919304590945695,   2.48862221210784) {};
\node[ele] (j1) at (1.64512903092519,    2.01677357394698) {};
\node[ele] (k1) at (1.24383442716840,    3.78580579686650) {};
\node[ele] (l1) at (1.58646607060451,    2.93089871716034) {};
\node[ele] (m1) at (0.874518089942222,   2.13017381291781) {};
\node[ele] (n1) at (1.46214937616206,    4.79124963787566) {};
\node[ele] (o1) at (2.63312304723148,    0.206799557528458) {};

\node (epedge) at (3.48, 3.58594837328201) {};
\node (ep) at (3.21, 3.71) {$\varepsilon$};
\draw[->] (b.center) -- (epedge.center) ;

\end{scope}

\begin{scope}[shift={(10,0)}]
    \node[ele] (a) at (2.8209087137275, 1.18093727264829) {};
    \node[ele] (b) at (2.95833231728538,	3.58594837328201) {};
    \node[ele] (c) at (0.311441915093690,	1.55228774215654) {};
    \node[ele] (d) at (1.69060866126193,	3.54516924716448) {};
    \node[ele] (e) at (1.31799287927655,	0.304557156170160) {};
    \node[ele] (f) at (1.61347687509477,	4.79163906483583) {};
    \node[ele] (g) at (2.27043233468693,	2.38135190156370) {};
    \node[ele] (h) at (1.20381702801476,	1.12542610389225) {};
    \node[ele] (i) at (0.919304590945695,	2.48862221210784) {};
    \node[ele] (j) at (1.64512903092519,	2.01677357394698) {};
    \node[ele] (k) at (1.24383442716840,	3.78580579686650) {};
    \node[ele] (l) at (1.58646607060451,	2.93089871716034) {};
    \node[ele] (m) at (0.874518089942222,	2.13017381291781) {};
    \node[ele] (n) at (1.46214937616206,	4.79124963787566) {};
    \node[ele] (o) at (2.63312304723148,	0.206799557528458) {};
    \foreach \firstnode in {a, b, c, d, e, f, g, h, i, j, k, l, m, n, o}{%
   \foreach \secondnode in {a, b, c, d, e, f, g, h, i, j, k, l, m, n, o}{%
   \setveclength{\mydist}{\firstnode}{\secondnode}
   \pgfmathparse{\mydist < 30 ? int(1) : int(0)}
   \ifnum\pgfmathresult=1
     \draw[red] (\firstnode.center) -- (\secondnode.center);
   \fi
  }
    \node[ele] (a) at (2.8209087137275, 1.18093727264829) {};
    \node[ele] (b) at (2.95833231728538,	3.58594837328201) {};
    \node[ele] (c) at (0.311441915093690,	1.55228774215654) {};
    \node[ele] (d) at (1.69060866126193,	3.54516924716448) {};
    \node[ele] (e) at (1.31799287927655,	0.304557156170160) {};
    \node[ele] (f) at (1.61347687509477,	4.79163906483583) {};
    \node[ele] (g) at (2.27043233468693,	2.38135190156370) {};
    \node[ele] (h) at (1.20381702801476,	1.12542610389225) {};
    \node[ele] (i) at (0.919304590945695,	2.48862221210784) {};
    \node[ele] (j) at (1.64512903092519,	2.01677357394698) {};
    \node[ele] (k) at (1.24383442716840,	3.78580579686650) {};
    \node[ele] (l) at (1.58646607060451,	2.93089871716034) {};
    \node[ele] (m) at (0.874518089942222,	2.13017381291781) {};
    \node[ele] (n) at (1.46214937616206,	4.79124963787566) {};
    \node[ele] (o) at (2.63312304723148,	0.206799557528458) {};
}
\begin{pgfonlayer}{background}
\path[fill=blue!30, draw] (g.center) to (j.center) to (l.center) to (g.center); 
\path[fill=blue!30, draw] (d.center) to (k.center) to (l.center) to (d.center); 
\path[fill=blue!30, draw] (m.center) to (i.center) to (j.center) to (m.center); 
\path[fill=blue!30, draw] (l.center) to (j.center) to (i.center) to (l.center); 
\end{pgfonlayer}

\end{scope}

\end{tikzpicture}
    \caption{The process of constructing a VR complex.}
    \label{fig: Rips}
\end{figure}

The first step in the process of TDA is to construct a \textit{simplicial complex}. Simplicial complexes are used as a way of attributing quantifiable shape to the data; they can be thought of as higher-dimensional analogues of graphs. In TDA, the vertices of the simplicial complexes are the observed data points, as seen in Fig. \ref{fig: Rips}. Simplicial complexes can be analysed to output the \textit{persistent homology}, a key \textit{topological invariant} that can be used to describe the structure of the data. The persistent homology can then be used to compare between different data sets, by quantifying their topological structures \cite{genomics}. A simplicial complex is made up of fundamental building blocks, called \textit{simplices}. The first four simplices are shown in Fig. \ref{fig:simplices}. Each vertex in the simplex is fully connected to all the other vertices, and the space enclosed by the vertices is part of that simplex. For instance, $\Delta^2$ encloses a two-dimensional area, $\Delta^3$ encloses a three-dimensional volume; this can be generalised for $\Delta^k$ enclosing a $k-$dimensional space between $(k+1)$ fully connected vertices.

There are many ways to construct a simplicial complex from point data. For simplification, only one method will be discussed within this paper, the \textit{Vietoris-Rips} (VR) complex \cite{carlsson2006algebraic}. The VR complex is constructed from point data, to output a corresponding simplicial complex, which can then be analysed. For the VR complex, $\text{VR}_\varepsilon$, let $(X,\partial_X)$ be a finite dimensional metric space, where $X$ denotes the set and $\partial_{X}$ is a metric on $X$. Then let $\varepsilon > 0$ be a fixed value, then \cite{chambers2010vietoris}:
\begin{enumerate}
    \item The vertices, $v \in X$, form the vertices in $\text{VR}_\varepsilon(X,\partial_X)$.
    \item A $k-$simplex is formed when $\partial_X(v_i,v_j) \leq 2\varepsilon, \ \forall i,j \leq k$ for some $\varepsilon > 0$.
\end{enumerate}
 
The process of constructing a VR complex is depicted in Fig. \ref{fig: Rips}, for some randomly-sampled data, and an arbitrary value of $\varepsilon$. The existence of a simplex is determined by how the balls intersect between vertices. For a VR complex, a simplex between some vertices is formed if the Euclidean distance between the all the vertices is less than $\varepsilon$.

From the simplicial complexes, the \textit{homology groups}, $H_k(X)$, can be determined. The homology groups are \textit{invariants} for the data set, $X$, where $k$ refers to the relevant dimension. Generally, the $k^{\text{th}}$ homology group encodes information about the number of $k-$dimensional holes in the data \cite{maclane2012homology, nash1988topology}. Under the rules of topology, discontinuities (voids) cannot be created or destroyed under continuous transformations \cite{intrototopologybertmendelsen, ghrist2018homological}. Therefore, the homology can be used to categorise and compare between simplicial complexes, and by extension, data sets \cite{boissonnat2018geometric, schutz1980geometrical, ghrist2014EAT}. From the homology, the \textit{Betti numbers} are defined as the \textit{rank} of the homology groups. If the Betti numbers for two topological spaces are different, these spaces are not topologically identical. If two spaces are not topologically similar, a \textit{continuous bijective map} between the spaces does not exist.
The zeroth Betti number, $\beta_0$, is the rank of the zeroth homology group \cite{nash1988topology}, $H_0(X)$, and refers to the number of connected sets in $X$. The first Betti number, $\beta_1$, is the rank of the first homology group, $H_1(X)$, and refers to the number of non-contractible holes present in $X$. The second Betti number, $\beta_2$, refers to the number of enclosed volumes in the topological space. This analogy carries on further for higher dimensions.

This description now raises the question: which length scale $\varepsilon$ is representative of the topology of the data? When constructing the VR complexes for the same data set, different values of epsilon, will result in different values for the Betti numbers. The hyper-parameter $\varepsilon$ determines the Betti numbers for that specific instance of some point cloud data. Additionally, when the feature present within the data is at a length scale less than $\varepsilon$ this feature will not be expressed, as $\varepsilon$ will span the feature. A problem arises here, as usually the feature scale is not known prior to analysis, and a manifold may have many multi-scale features. The answer to this problem, is to vary $\varepsilon$ and see how the Betti numbers evolve and \textit{persist}. Obtaining the homology for a single value of $\varepsilon$ provides very limited information, because of potentially-varying feature length scales in the manifold. For this reason, it is vital to consider how homological features persist as $\varepsilon$ is varied. This process of varying $\varepsilon$ does not bias any disk size, as all are being considered. This process will give an initial value, $\varepsilon_{\text{min}}$, where a specific homological feature comes to life and $\varepsilon_{\text{max}}$, where the feature is no longer considered for that simplicial complex. This range of values $[\varepsilon_{\text{min}},\varepsilon_{\text{max}}]$ is called the \textit{persistent interval} for that homological feature. Each persistent interval is attributed a Betti number. The set of persistent intervals, with their accompanying Betti numbers forms an object called the \textit{persistent homology} (PH). PH is invariant for a data set and contains enough information to represent the topology of a point cloud.

The space of persistent homologies forms a \textit{metric space}; the distance between the PHs is a measure of similarity of two PHs. As the persistent intervals are invariant for a manifold, the data manifolds can be represented by their persistent homology. This notion of a metric space allows one to compare the similarity of manifolds. Metrics between barcodes are well established and the one used in this report is the \textit{$p-$Wasserstein distance}, $\partial_{W_p}$. 
$$\partial_{W_p}(B_1,B_2) = \left(\text{inf} \sum_{Z\in B_1}d_\infty(Z, \phi(Z))^p \right)^{\frac{1}{p}}$$
Where $B_1$ and $B_2$ are PH data, $p>0$ is a weighting, $\phi$ is a matching between $B_1$ and $B_2$, $Z$ is a persistent interval in $B_1$, and $d_\infty$ is the supremum metric \cite{genomics}.

\subsection{Time-Delay Embedding}
Given a time-varying series, $f:t \rightarrow \mathbb{R}$, the time-delay embedding can be stacked $d$ times, each with a delay $\alpha$, to give a new embedding $\phi : t \rightarrow \mathbb{R}^d$ where the new embedding is represented by,

\begin{equation}
    \phi(t) = (f(t), f(t+\alpha), ..., f(t+(d-1)\alpha)).
\end{equation}

From this embedding, the time series now has a topology induced. This embedding is a reconstruction of the topology from a one-dimensional time series into any desired dimension. The first application of a time delay embedding was to give geometry to a time series \cite{packard1980geometry}. Since then, it has been shown that time-delay embeddings can reconstruct the topology of a dynamic attractor, under certain circumstances \cite{takens1981detecting}.

\begin{figure}[H]
    \centering
    \includegraphics[width = 0.5\textwidth]{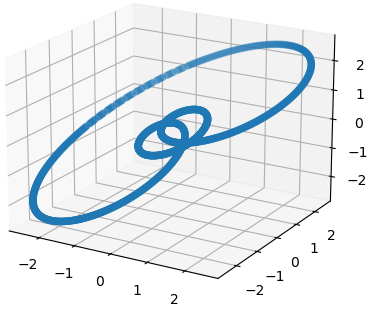}
    \caption{Data generated from the time series $y(t) = \sin(t) + \sin(2t) + \sin(3t)$, embedded in $3D$ at $\alpha=5$.}
    \label{fig:tdetoy}
\end{figure}

The delay embedding highlights periodic recurrent features in the time series.  Recurrent behaviour will be highlighted in the time-delay embedding as a loop. Persistent homology can then be used to quantify the size and number of these loops. Fig. \ref{fig:tdetoy} shows time-delay embedding with multiple frequency sines, different sized loops are formed as a result of these different frequencies.

If the delay is too small, there will not be sufficient information to form meaningful topology, and the reconstruction will be too similar to a straight diagonal line, as the local change is too small. Conversely, using a too-large delay will result in nonsense, as the gap between the readings will be too large and will not show a local change over the manifold.

\section{Results}
\subsection{Linear and Nonlinear Cointegration of $\omega_2$}
The effective procedure of performing the TDA is to embed all the time-series data into a common vector space. This idea assumes that if the time series all have consistent topology, there should be no outliers when compared to one another. This assumption can be used to compare the topologies of the embedded time series.

\begin{figure}[H]
    \centering
    \subfloat[{Non-cointegrated data.}]{\includegraphics[width = 0.47\textwidth]{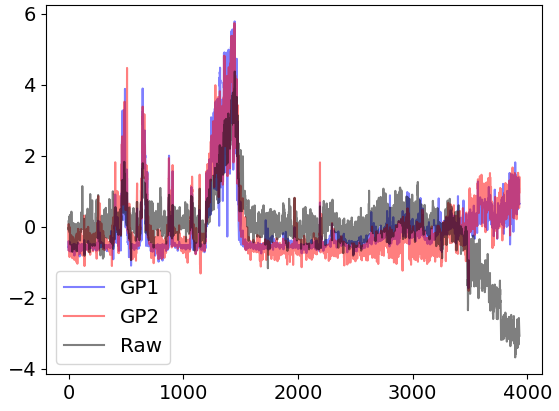}}
    \qquad
    \subfloat[{Cointegrated data.}]{\includegraphics[width = 0.47\textwidth]{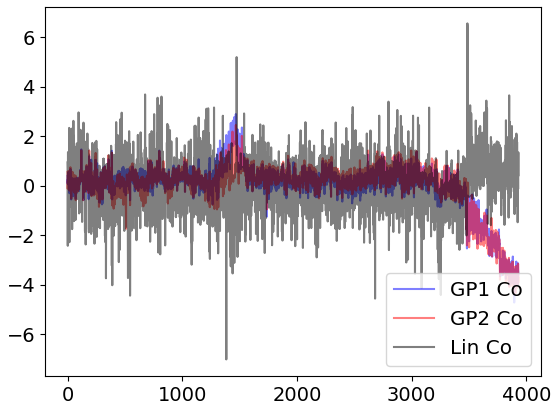}}
    \qquad
    \caption{Univariate plots of all the time series representing $\omega_2$.}
    \label{fig: timeseries}
\end{figure}

There are six distinct time series here, all representing $\omega_2$, that are going to be analysed. All of the time series have been standardised prior to embedding. These time series can be seen in Fig. \ref{fig: timeseries}. The significance of each time series, is: 
\begin{enumerate}[GP1 CO]
    \item[\textbf{RAW:}] This is the raw time series for $\omega_2$, calculated from modal analysis \cite{peeters2001one}, with NaNs removed.
    \item[\textbf{GP1:}] From the \textbf{RAW} time series, a GP was fitted using $\omega_1, \omega_3, \omega_4$ and was trained over the first 1000 data points, by index.
    \item[\textbf{GP2:}] From the \textbf{RAW} time series, a GP was fitted using $\omega_1, \omega_3, \omega_4$ and was trained over a 1000 points that included a large temperature drop.
    \item[\textbf{LIN CO:}] This is the most stationary residual, when the first four natural frequencies from the Z24 data set are cointegrated.
    \item[\textbf{GP1 CO:}] This is the residual series from \textbf{GP1}, $\epsilon_1 = \omega_2 - \textbf{GP1}(\omega_1, \omega_3, \omega_4)$.
    \item[\textbf{GP2 CO:}] This is the residual series from \textbf{GP2}, $\epsilon_1 = \omega_2 - \textbf{GP2}(\omega_1, \omega_3, \omega_4)$.
\end{enumerate}

\begin{figure}
    \centering
    \includegraphics[width=0.6\textwidth]{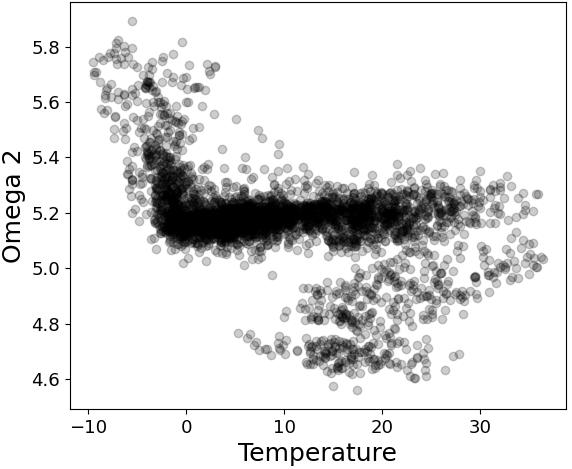}
    \caption{The second natural frequency vs temperature, showing a nonlinear relationship.}
    \label{fig:nonlinomega}
\end{figure}

The Z24 Bridge data set is well studied in the SHM community, and the second natural frequency is known to be nonlinear with respect to temperature \cite{shi2019cointegration, cross2011approaches}, this can be seen in Fig. \ref{fig:nonlinomega}. The reasoning to include both nonlinear and linear cointegration schemes on the second natural frequency is to assess the topological changes in the nonlinear time series under both implementations of cointegration.

All the time series listed above were embedded into three-dimensional space. TDA is used to calculate the persistent homology of the embedded time series. The persistent homologies can then be compared to one another by use of the Wasserstein distance; the greater the Wasserstein distance, the more dissimilar the topologies are.

\begin{figure}[H]
    \centering
    \includegraphics[width=\textwidth]{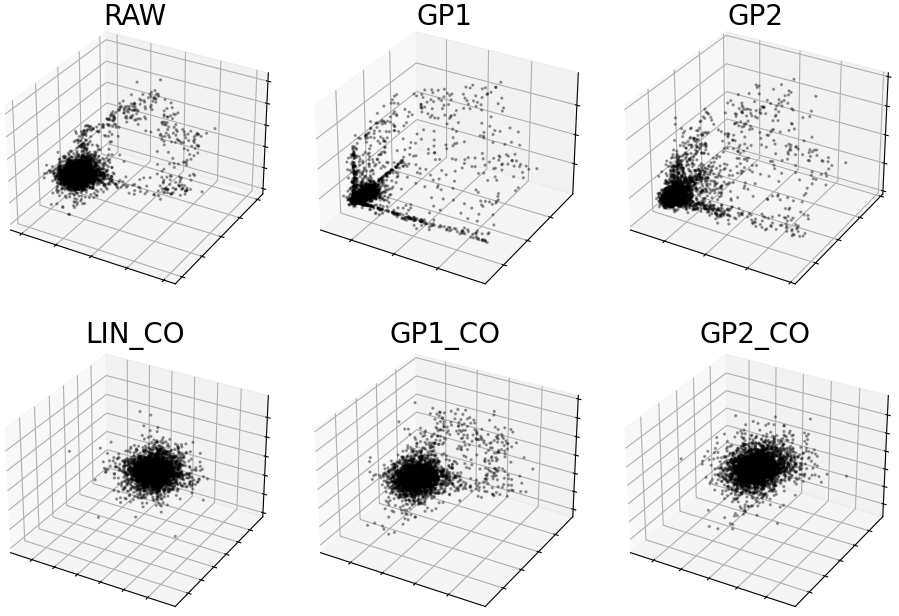}
    \caption{The time series data embedded into $3D$ at a delay $\alpha=75$.}
    \label{fig:3dembedded}
\end{figure}

As can be seen in Fig. \ref{fig:3dembedded}, there are two prominent topologies appearing; a torus for the non-cointegrated data, and an open-ball for the cointegrated data. The torus arises as a result of the cyclic temperature effects. This topology is most clear in the plot for RAW. The torus is still present in both GP1 and GP2, but is not as pronounced as the GPs are trained on the natural frequencies $\omega_1, \omega_3, \omega_4$. Therefore, nonlinear temperature effects are not represented.

The open-ball topology can be seen for the cointegrated data in Fig. \ref{fig:3dembedded}. This topology is formed because the residuals are white-noise processes, therefore the embedding results in a Gaussian cluster. For the residuals from GP1, it can be seen that there are temperature effects still visible in the residuals. This is a direct result of GP1 being trained over the the first 1000 points in the data set, where there is no large temperature changes present. Therefore, GP1 does not adequately purge the nonlinear temperature effects from the residual series. This is different to GP2, which is trained over the large temperature drop, and there is no torus present in the embedded data, indicating that this is a more appropriate GP to detrend the residual series. 

\begin{table}
    \centering
    \begin{tabular}{r|ccc|ccc}
     & \textbf{RAW} & \textbf{GP1} & \textbf{GP2} & \textbf{LIN CO} & \textbf{GP1 CO} & \textbf{GP2 CO} \\
     \hline
    \textbf{RAW}     &   0.00 & 131.00 &  80.74 &  196.99 &  103.05 &  174.39 \\
    \textbf{GP1}     & 131.00 &   0.00 &  63.32 &  272.77 &  194.04 &  253.30 \\
    \textbf{GP2}     & 80.74  &  63.32 &   0.00 &  230.11 &  143.88 &  209.59 \\
    \hline
    \textbf{LIN CO}  & 197.00 & 272.79 & 230.12 &    0.00 &  106.30 &   36.89 \\
    \textbf{GP1 CO}  & 103.05 & 194.04 & 143.88 &  106.31 &    0.00 &   84.63 \\
    \textbf{GP2 CO}  & 174.39 & 253.30 & 209.59 &   36.89 &   84.63 &    0.00 \\
    \end{tabular}
    \caption{The Wasserstein distances when two data sets are compared, at a delay $\alpha=75$ in an embedding dimension of 3.}
    \label{tab:wassdelay_75}
\end{table}

Table \ref{tab:wassdelay_75} shows the Wasserstein distances between all of the embedded time-series. Table \ref{tab:wassdelay_75} indicates that GP2 is a better fit to RAW than GP1. This is likely due to the poor representation of GP1 over the large temperature drop. This fit is also shown in the nonlinear cointegration comparisons, as now the embedding for GP1 CO is more similar to RAW than GP2 CO, as GP1 has detrended by a lesser extent.   

If Table \ref{tab:wassdelay_75} is considered as a 2x2 block matrix of 3x3 matrices (indicated by the lines partitioning the table). The top-left 3x3 block shows comparisons between the three toroidal topologies. The off-diagonal blocks show comparisons between the tori and the balls. The bottom-right block shows comparisons between the balls. Therefore, by this reasoning it is expected that data in the leading diagonal blocks are relatively small, and the the off diagonal blocks relatively large. This schema is shown in Table \ref{tab:tabstruct}.

\begin{table}
    \centering
    \begin{tabular}{c|c|c}
            & \textbf{Tori} \ & \textbf{Balls}  \\
            \hline
            \textbf{Tori} \  & Low \ & High \\
            \hline
            \textbf{Balls} \ & High \ & Low
    \end{tabular}
    \caption{Overview of the structure of Table \ref{tab:wassdelay_75}.}
    \label{tab:tabstruct}
\end{table}

\subsection{Linear Cointegration}
\begin{figure}
    \centering
    \includegraphics[width = 0.6\textwidth]{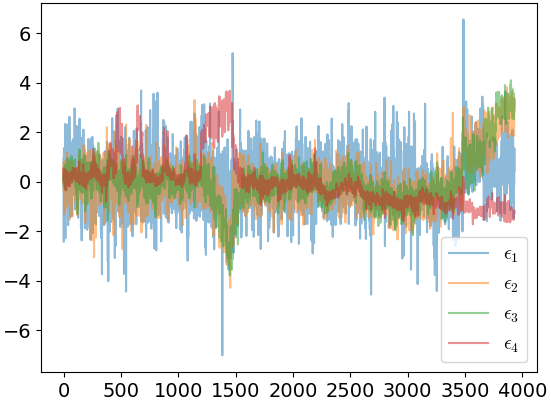}
    \caption{The linear residual series calculated from the first four natural frequencies.}
    \label{fig:lincoinres}
\end{figure}

The same process is now applied to the linear cointegrated residual series, calculated from  the first four natural frequencies. The linear cointegrated residuals are shown in Fig. \ref{fig:lincoinres}. Again, all the series have been standardised before embedding.

The  Wasserstein distances between the embedded linearly cointegrated residual series can be found in Table \ref{tab:fulllincoin}. For these results, a delay of $75$, and an embedding dimension of three were used. 

\begin{table}
    \centering
    \begin{tabular}{r|cccc|cccc}
     & $\omega_1$ & $\omega_2$ & $\omega_3$ & $\omega_4$ & $\epsilon_1$ & $\epsilon_2$ & $\epsilon_3$ & $\epsilon_4$\\
     \hline
    $\omega_1$     &  0.00 &  15.60 &  61.28 &  85.40 &   419.97 &   325.38 &   258.48 &   193.53 \\
    $\omega_2$     & 15.60 &   0.00 &  47.53 &  72.82 &   421.11 &   323.32 &   254.47 &   186.25 \\
    $\omega_3$     & 61.28 &  47.53 &   0.00 &  29.24 &   407.34 &   301.73 &   227.74 &   152.38 \\
    $\omega_4$ & 85.40 &  72.82 &  29.24 &   0.00 &   395.48 &   285.63 &   208.31 &   129.74 \\
    \hline
    $\epsilon_1$ & 419.97 & 421.11 & 407.34 & 395.49 &     0.00 &   128.08 &   211.80 &   286.00 \\
    $\epsilon_2$ & 325.38 & 323.32 & 301.73 & 285.63 &   128.08 &     0.00 &    89.48 &   170.15 \\
    $\epsilon_3$ & 258.48 & 254.47 & 227.74 & 208.31 &   211.80 &    89.47 &     0.00 &    88.83 \\
    $\epsilon_4$ & 193.54 & 186.25 & 152.39 & 129.75 &   286.00 &   170.15 &    88.84 &     0.00 \\
    \end{tabular}
    \caption{Wasserstein distances for the embedded time series before and after cointegration, for a delay of 75 and an embedding dimension of 3.}
    \label{tab:fulllincoin}
\end{table}
The pattern, shown in Table \ref{tab:tabstruct} is still visible in the linear cointegrated data. Inside the two off diagonal blocks, as progressive residual series are used, the Wasserstein distances are decreasing. This is to be expected as the first residual series is most stationary. Therefore, more of the topology has been purged in the first residual series when compared to later ones. As successive residual series are considered, some nonlinear trends are breaking through. As a result, the topology is becoming closer to the raw data embedding.

\section{Conclusion}
With respect to the Z24 data set, topological methods have been used to analyse time series data, before and after cointegration. When temperature effects are not adequately represented in the GP model, the topology of the original time series seems to break through into the residuals. This work has shown that time-delay embeddings and TDA are able to quantify how well cointegration has detrended a time series, by analysing the topology of the residuals.

Further work on cointegration and TDA will generalise the ideas found within this paper to more general data. Data will be generated that has a more complex topology when a time-delay embedding is computed. A nonlinear cointegration scheme will be applied and assessed with TDA.

\section*{Acknowledgements}
The authors would like to thank the UK EPSRC for funding via the Established Career Fellowship EP/R003645/1 and the Programme Grant EP/R006768/1.



\bibliographystyle{ieeetr}
\bibliography{EWSHM_319.bib}

\end{document}